\documentstyle[amscd,amssymb,xypic,12pt]{amsart}
\xyoption{all}

\newcommand{\nc}{\newcommand}

\nc{\exto}[1]{\stackrel{#1}{\longrightarrow}}
\nc{\lan}{\big\langle}
\nc{\ran}{\big\rangle}

\nc{\C}{{\mathbb{C}}}
\nc{\HH}{{\mathbb{H}}}
\nc{\PP}{{\mathbb{P}}}
\nc{\QQ}{{\mathbb{Q}}}
\nc{\ZZ}{{\mathbb{Z}}}

\nc{\CA}{{\mathcal{A}}}
\nc{\CC}{{\mathcal{C}}}
\nc{\D}{{\mathcal{D}}}
\nc{\CE}{{\mathcal{E}}}
\nc{\CF}{{\mathcal{F}}}
\nc{\CH}{{\mathcal{H}}}
\nc{\CL}{{\mathcal{L}}}
\nc{\CM}{{\mathcal{M}}}
\nc{\CO}{{\mathcal{O}}}
\nc{\CQ}{{\mathcal{Q}}}
\nc{\CU}{{\mathcal{U}}}
\nc{\CV}{{\mathcal{V}}}

\nc{\FE}{{\mathfrak{E}}}
\nc{\FL}{{\mathfrak{L}}}
\nc{\FM}{{\mathfrak{M}}}
\nc{\FS}{{\mathsf{S}}}

\nc{\sfc}{{\mathsf{c}}}
\nc{\sfch}{{\mathsf{ch}}}

\nc{\SK}{{\mathsf{K}}}
\nc{\SO}{{\mathsf{O}}}
\nc{\SPV}{{\mathsf{S}^+\mathsf{V}}}
\nc{\SMV}{{\mathsf{S}^-\mathsf{V}}}
\nc{\SPMV}{{\mathsf{S}^\pm\mathsf{V}}}
\nc{\SX}{{S_X}}
\nc{\SY}{{S_Y}}
\nc{\phipsi}{{q}}

\nc{\pim}{{\pi_-}}
\nc{\pip}{{\pi_+}}

\nc{\TE}{{\tilde{\CE}}}
\nc{\TW}{{\tilde{W}}}
\nc{\txi}{{\tilde{\xi}}}
\nc{\tp}{{\tilde{p}}}
\nc{\tq}{{\tilde{q}}}

\nc{\HX}{{\hat{X}}}
\nc{\hxi}{{\hat{\xi}}}

\nc{\UH}{{\mathcal{H}}}

\nc{\TM}{{\widetilde{M}}}
\nc{\TX}{{\widetilde{X}}}
\nc{\TY}{{\widetilde{Y}}}
\nc{\barf}{{\bar{f}}}
\nc{\tf}{{\tilde{f}}}
\nc{\hf}{{\hat{f}}}

\nc{\LRA}{\Leftrightarrow}
\nc{\RA}{\Rightarrow}
\nc{\lotimes}{\mathbin{\mathop{\otimes}\limits^{\mathbb{L}}}}
\nc{\CExt}{\mathop{\mathcal{E}\mathit{xt}}\nolimits}
\nc{\CHom}{\mathop{\mathcal{H}\mathit{om}}\nolimits}
\nc{\RH}{\mathop{{\mathsf{R}}\Gamma}\nolimits}
\nc{\RHom}{\mathop{\mathsf{RHom}}\nolimits}
\nc{\RCHom}{\mathop{\mathsf{R}\mathcal{H}\mathit{om}}\nolimits}
\nc{\RG}{\mathop{\mathsf{R\Gamma}}\nolimits}
\nc{\Hom}{\mathop{\mathsf{Hom}}\nolimits}
\nc{\Ext}{\mathop{\mathsf{Ext}}\nolimits}
\nc{\End}{\mathop{\mathsf{End}}\nolimits}
\nc{\Tor}{\mathop{\mathsf{Tor}}\nolimits}
\nc{\Hilb}{\mathop{\mathsf{Hilb}}\nolimits}
\nc{\Spec}{\mathop{\mathsf{Spec}}\nolimits}
\nc{\Pic}{\mathop{\mathsf{Pic}}\nolimits}

\nc{\Tw}{\mathop{\mathsf{Tw}}\nolimits}
\nc{\Ker}{\mathop{\mathsf{Ker}}\nolimits}
\nc{\Coker}{\mathop{\mathsf{Coker}}\nolimits}
\nc{\codim}{\mathop{\mathsf{codim}}\nolimits}
\nc{\sing}{{\mathsf{sing}}}
\nc{\supp}{{\mathsf{supp}}}
\nc{\rank}{{\mathsf{rank}}}
\nc{\Pf}{{\mathsf{Pf}}}
\nc{\Gr}{{\mathsf{Gr}}}
\nc{\LGr}{{\mathsf{LGr}}}
\nc{\Fl}{{\mathsf{Fl}}}
\nc{\Bl}{{\mathsf{Bl}}}
\nc{\GL}{{\mathsf{GL}}}
\nc{\Spin}{{\mathsf{Spin}}}
\nc{\ev}{{\mathsf{ev}}}
\nc{\id}{{\mathsf{id}}}

\nc{\Cubics}{{{\mathcal{S}}_3}}
\nc{\VFT}{{{\mathcal{S}}_{14}}}
\nc{\VFTE}{{{\mathcal{N}}_{\mathrm{reg,sm}}}}
\nc{\MX}{{\CM_X}}
\nc{\MY}{{\CM_Y}}
\nc{\MYE}{{\CM_{Y,\CE}}}

\theoremstyle{plain}

\newtheorem{theorem}{Theorem}[section]
\newtheorem{lemma}[theorem]{Lemma}
\newtheorem{proposition}[theorem]{Proposition}
\newtheorem{corollary}[theorem]{Corollary}

\theoremstyle{definition}

\theoremstyle{remark}

\newtheorem{remark}[theorem]{Remark}

\newenvironment{proof}{\noindent{\sf Proof:}}{\qed\medskip}

\title{Derived category of $V_{12}$ Fano threefolds}
\author{Alexander Kuznetsov}
\address{
Algebra Section, Steklov Mathematical Institute,
Russian Academy of Sciences,
 8 Gubkin str., Moscow 119991, Russia
}
\email{akuznet@@mi.ras.ru, sasha@@kuznetsov.mccme.ru}
\date{}

\begin{document}

\maketitle

\section{Introduction}

A $V_{12}$ Fano threefold is a smooth Fano threefold $X$ of index $1$ 
with $\Pic X = \ZZ$ and $(-K_X)^3=12$, see \cite{Is1,IP}. 
Let $X$ be a $V_{12}$ threefold. It was shown by Mukai \cite{Mu1} 
that $X$ admits an embedding into a connected component $\LGr_+(V)$
of the Lagrangian Grassmannian $\LGr(V)$ of Lagrangian ($5$-dimensional) 
subspaces in a vector space $V = \C^{10}$ with respect to a nondegenerate
quadratic form~$Q$, and moreover, $X = \LGr_+(V)\cap\PP^8$. Let $\CO_X$ be
the structure sheaf and let $\CU_+$ denote the restriction to $X$ of 
the tautological ($5$-dimensional) subbundle from $\LGr_+(V) \subset \Gr(5,V)$. 
Then it is easy to show that $(\CU_+,\CO_X)$ is an exceptional pair in 
the bounded derived category of coherent sheaves on $X$, $\D^b(X)$. 
Therefore, triangulated subcategory $\langle\CU_+,\CO_X\rangle$ 
generated by the pair is admissible and there exists a semiorthogonal 
decomposition $\D^b(X) = \langle\CU_+,\CO_X,\CA_X\rangle$, where 
$\CA_X = {}^\perp\langle\CU_+,\CO_X\rangle \subset \D^b(X)$ 
is the orthogonal subcategory. The main result of this note is 
an equivalence $\CA_X\cong\D^b(C^\vee)$, where $C^\vee$ is 
a curve of genus~$7$. 

The curve $C^\vee$ arising in this way in fact is nothing
but the orthogonal section of the Lagrangian Grassmannian considered 
by Iliev and Markushevich \cite{IM}. Recall that the components $\LGr_+(V)$
and $\LGr_-(V)$ of the Lagrangian Grassmannian $\LGr(V)$ lie in the dual
projective spaces $\PP(\SPV)$ and $\PP(\SMV)$ respectively, where $\SPMV$
are the spinor ($16$-dimensional) representations of the corresponding 
spinor group $\Spin(V)$. So, with any linear subspace $\PP^8\subset\PP(\SPV)$ 
one can associate its orthogonal $(\PP^8)^\perp = \PP^6\subset\PP(\SMV)$ and 
consider following~\cite{IM} the orthogonal section 
$C^\vee := \LGr_-(V)\cap\PP^6$, which can be shown to be 
a smooth genus~$7$ curve, whenever $X$ is smooth.

Further, Iliev and Markushevich explained in \cite{IM} the intrinsic meaning 
of the curve~$C^\vee$ associated to the threefold~$X$. 
They have shown that it is isomorphic to the moduli space of
stable rank $2$ vector bundles on $X$ with $c_1 = 1$, $c_2 = 5$.
Considering a universal bundle ${\CE_1}$ on $X\times C^\vee$ we obtain
the corresponding kernel functor $\Phi_{\CE_1}:\D^b(C^\vee) \to \D^b(X)$.
It follows from \cite{BO} and \cite{IM} that $\Phi_{\CE_1}$ is fully faithful.
Moreover, it can be shown that its image is contained in the orthogonal
subcategory $\CA_X = {}^\perp\langle\CU_+,\CO_X\rangle \subset \D^b(X)$.
Thus, it remains to check that $\Phi_{\CE_1}:\D^b(C^\vee)\to\CA_X$ is 
essentially surjective. 

To prove the surjectivity of the functor $\Phi_{\CE_1}$ we use 
the following approach. Take arbitrary smooth hyperplane section 
$X\supset S:=X \cap\PP^7=\LGr_+(V)\cap\PP^7$ and consider the orthogonal 
section $S^\vee = \LGr_-(V)\cap(\PP^7)^\perp$. Then both $S$ and $S^\vee$ 
are $K3$ surfaces, moreover $S^\vee$ is smooth and $C^\vee$ is a hyperplane 
section of $S^\vee$. Iliev and Markushevich have shown in \cite{IM} that 
the moduli space of stable rank $2$ vector bundles on $S$ with $c_1 = 1$, 
$c_2 = 5$ is isomorphic to $S^\vee$, so we can again consider a universal 
bundle $\CE_2$ on $S\times S^\vee$ and the corresponding kernel functor 
$\Phi_{\CE_2}:\D^b(S^\vee) \to \D^b(S)$. Again it follows from \cite{BO} 
and \cite{IM} that $\Phi_{\CE_2}$ is fully faithful, hence an equivalence 
by \cite{Br1}. Further, it is clear that we have an isomorphism
${\CE_1}_{|S\times C^\vee} \cong {\CE_2}_{|S\times C^\vee}$, hence 
the composition of $\Phi_{\CE_1}$ with pushforward from $C^\vee$ to $S^\vee$ 
coincides with the composition of $\Phi_{\CE_2}$ with restriction 
from $X$ to $S$: $\alpha^*\circ\Phi_{\CE_1} \cong \Phi_{\CE_2}\circ\beta_*$,
where $\alpha:S\to X$ and $\beta:C^\vee\to S^\vee$ are the embeddings.
The crucial observation however is that the bundles $\CE_1$ on $X\times C^\vee$
and $\CE_2$ on $S\times S^\vee$ can be glued on $X\times S^\vee$ so that 
the kernel functor $\D^b(X) \to \D^b(S^\vee)$ corresponding to the glueing 
vanishes on the subcategory $\CA_X$. In other words, we have
$(\beta_*\circ\Phi_{\CE_1}^*)_{|\CA_X} \cong 
(\Phi_{\CE_2}^*\circ\alpha^*)_{|\CA_X}$,
where $\Phi_{\CE_1}^*$ and $\Phi_{\CE_2}^*$ are the left adjoint functors.
%
%
Now the proof goes as follows. Take an object $F\in\CA_X$,
orthogonal to the image of $\Phi_{\CE_1}$. Then $\Phi_{\CE_1}^*(F)=0$.
Hence $(\Phi_{\CE_2}^*\circ\alpha^*)(F)=(\beta_*\circ\Phi_{\CE_1}^*)(F)=0$. 
But $\Phi_{\CE_2}$ is an equivalence, hence $\Phi_{\CE_2}^*$ is 
an equivalence, hence $\alpha^*(F)=0$. Since these arguments
apply to {\em any}\/ smooth hyperplane section $S\subset X$,
it follows that the restriction of such $F$ to any smooth hyperplane
section is zero, but this immediately implies that $F=0$.

Having in mind the semiorthogonal decomposition 
$\D^b(X) = \langle\CU_+,\CO_X,\D^b(C^\vee)\rangle$ 
one can informally say that the nontrivial part of the derived
category $\D^b(X)$ is described by the curve $C^\vee$. Therefore,
the curve $C^\vee$ should appear in all geometrical questions
related to $X$. As a demonstration of this phenomenon we show
that the Fano surface of conics on $X$ is isomorphic to the
symmetric square of $C^\vee$. This fact was known to Iliev and 
Markushevich, see \cite{IM2}, however we decided to include our proof 
into the paper for two reasons: it demonstrates very well how the above 
semiorthogonal decomposition can be used, and, moreover, the same approach
allows to investigate any other moduli space on $X$.

The paper is organised as follows. In section 2 we recall briefly 
results of \cite{IM}. In section 3 we give an explicit description
of universal bundles on $X\times C^\vee$ and $S\times S^\vee$ and
of their glueing on $X\times S^\vee$. In section 4 we give necessary
cohomological computations. In section 5 we consider the derived 
categories and prove the equivalence $\CA_X \cong \D^b(C^\vee)$.
Finally, in section~6 we investigate conics on $X$ and prove
that the Fano surface $F_X$ is isomorphic to $S^2C^\vee$.

\bigskip

{\sc Acknowledgements.}
I am grateful to Atanas Iliev and Dmitry Markushevich
for valuable communications on $V_{12}$ threefolds and genus~$7$ curves 
and to Dmitry Orlov and Alexei Bondal for useful discussions.
I was partially supported by RFFI grants 02-01-00468 and 02-01-01041 
and INTAS-OPEN-2000-269. The research described in this work 
was made possible in part by CRDF Award No. RM1-2405-MO-02.

\section{Preliminaries}

Fix a vector space $V=\C^{10}$ and a quadratic nondegenerate form $Q$ on $V$. 
Let $\SPV$, $\SMV$ denote the spinor ($16$-dimensional) representations 
of the spinor group $\Spin(Q)$. Recall that the spaces $\SPMV$ coincide 
with the (duals of the) spaces of global sections of the ample generators 
of the Picard group of connected components $\LGr_\pm(V)$ of the Lagrangian 
Grassmanian of $V$ with respect to $Q$. In particular, we have canonical 
embeddings $\LGr_\pm(V) \to \PP(\SPMV)$.

Choose a pair of subspaces $A_8 \subset A_9 \subset \SPV$, 
$\dim A_i = i$, and consider the intersections
\begin{equation}\label{CSX}
\begin{array}{rrr}
S = \LGr_+(V) \cap \PP(A_8) \subset \PP(\SPV),\\
X = \LGr_+(V) \cap \PP(A_9) \subset \PP(\SPV).
\end{array}
\end{equation}
It is easy to see that if $X$ is smooth then $X$ is a $V_{12}$
threefold, and if $S$ is smooth then $S$ is a polarized $K3$ surface
of degree~$12$.

\begin{theorem}[\cite{Mu1}]
If $X$ is a $V_{12}$ Fano threefold, and $S \subset X$ is its smooth
$K3$ surface section, then there exists a pair of subspaces
$A_8 \subset A_9 \subset \SPV$, such that $S$ and $X$ are obtained 
by~$(\ref{CSX})$.
\end{theorem}

Recall that the spinor representations $\SMV$ and $\SPV$ are canonically 
dual to each other, and denote by $B_7 \subset B_8 \subset \SMV$ 
the orthogonal subspaces, 
$$
B_i = A_{16-i}^\perp \subset {\SPV}^* \cong \SMV,
$$
and consider the dual pair
\begin{equation}\label{CSX1}
\begin{array}{rrr}
C^{\vee} = \LGr_-(V) \cap \PP(B_7) \subset \PP(\SMV),\\
S^{\vee} = \LGr_-(V) \cap \PP(B_8) \subset \PP(\SMV).
\end{array}
\end{equation}
Again, it is easy to see that if $S^\vee$ is smooth then $S^\vee$ 
is a polarized $K3$ surface of degree~$12$, and if $C^\vee$ is smooth
then $C^\vee$ is a canonically embedded curve of genus~$7$.

We denote by $H_X$, $L_X$, and $P_X$ the classes of a hyperplane section,
of a line, and of a point in $H^\bullet(X,\ZZ)$. The same notation is used
for varieties $S$, $S^\vee$ and $C^\vee$. For example, 
$P_{S^\vee}\in H^4(S^\vee,\ZZ)$ stands for the class of a point on $S^\vee$.

Let $\CU_+$, $\CU_-$ denote the tautological subbundles on 
$\LGr_+(V) \subset \Gr(5,V)$, $\LGr_-(V) \subset \Gr(5,V)$
respectively, and by ${\CU_+}_x$, ${\CU_-}_y$ their fibers
at points $x\in\LGr_+(V)$, $y\in\LGr_-(V)$ respectively.

Recall the relation between the canonical duality of $\LGr_\pm(V)$
and the intersection of subspaces.

\begin{lemma}[\cite{IM}]\label{di}
Let $x\in\LGr_+(V) \subset \PP(\SMV)$, $y\in\LGr_-(V) \subset \PP(\SPV)$
and denote by $\langle-,-\rangle$ the duality pairing on $\SPV\times\SMV$.
Then 
$$
\begin{array}{lcr}
\langle x,y \rangle \ne 0 & \LRA & 
{\CU_+}_x \cap {\CU_-}_y = 0,
\smallskip\\
\langle x,y \rangle = 0 & \RA  & 
\dim \left({\CU_+}_x \cap {\CU_-}_y\right) \ge 2.
\end{array}
$$
\end{lemma}

It follows that for any $(x,y)\in X\times C^\vee$ or 
$(x,y) \in S\times S^\vee$ we have 
$\dim \left({\CU_+}_x \cap {\CU_-}_y\right) \ge 2$.

\begin{lemma}[\cite{IM}]\label{smoothness}
We have the following equivalences:

\noindent $(i)$
$C^\vee$ is smooth $\LRA$ $X$ is smooth $\LRA$ 
for all $x\in X$, $y \in C^\vee$ we have $\dim({\CU_+}_x \cap {\CU_-}_y) = 2$;

\noindent $(ii)$
$S^\vee$ is smooth $\LRA$ $S$ is smooth $\LRA$ 
for all $x\in S$, $y \in S^\vee$ we have $\dim({\CU_+}_x \cap {\CU_-}_y) = 2$.
\end{lemma}

The following theorem reveals the intrinsic meaning of the curve $C^\vee$ 
and of the surface $S^\vee$ in terms of $X$ and $S$ respectively.

\begin{theorem}[\cite{IM}]\label{im2}
$(i)$ The curve $C^\vee$ is the fine moduli space of stable rank $2$
vector bundles $E$ on $X$ with $\sfc_1(E)=H_X$, $\sfc_2(E)=5L_X$. 
If $E_y$, $E_{y'}$ are the bundles on $X$ corresponding 
to points $y,y'\in C^\vee$, then
$$
\Ext^p(E_y,E_{y'}) = \begin{cases}
\C, & \text{for $p=0,1$ and $y=y'$}\\
0, & \text{otherwise}
\end{cases}
$$

\noindent
$(ii)$ The surface $S^\vee$ is the fine moduli space of stable rank $2$
vector bundles $E$ on $S$ with $\sfc_1(E)=H_S$, $\sfc_2(E)=5P_S$.
If $E_y$, $E_{y'}$ are the bundles on $S$ corresponding 
to points $y,y'\in S^\vee$, then
$$
\Ext^p(E_y,E_{y'}) = \begin{cases}
\C, & \text{for $p=0,2$ and $y=y'$}\\
\C^2, & \text{for $p=1$ and $y=y'$}\\
0, & \text{otherwise}
\end{cases}
$$
\end{theorem}

\section{The universal bundles}

Consider one of the following two products
$$
\text{either}\quad W_1 = X \times C^\vee,
\qquad\qquad\text{or}\quad
W_2 = S \times S^\vee.
$$
Denote by $\CU_+$ and $\CU_-$ the pullbacks of the tautological subbundles 
on $\LGr_+(V)$ and $\LGr_-(V)$ to $W_i \subset \LGr_+(V)\times\LGr_-(V)$,
and consider the following natural composition of morphisms 
of vector bundles on $W_i$
$$
\xi_i: \CU_- \to V\otimes\CO_{W_i} \stackrel{Q}\cong
V^*\otimes\CO_{W_i} \to \CU_+^*.
$$

\begin{lemma}
If $X$ {\rm(}resp.\ $S${\rm)} is smooth then the rank of $\xi_1$ 
{\rm(}resp.\ $\xi_2${\rm)} equals $3$ at every point of $W_1$
{\rm(}resp.\ $W_2${\rm)}.
\end{lemma}
\begin{proof}
Since the kernel of the natural projection
$V^*\otimes\CO_{\LGr_+(V)} \to \CU_+^*$ 
equals $\CU_+$, it suffices to show that 
for all points $(x,y)\in W_i\subset\LGr_+(V)\times\LGr_-(V)$
we have $\dim({\CU_+}_x\cap{\CU_-}_y) = 2$ which follows
from lemma~\ref{smoothness}. 
\end{proof}

\begin{lemma}
We have $\Ker \xi_i \cong (\Coker \xi_i)^*$.
\end{lemma}
\begin{proof}
We have the following commutative diagram with exact rows:
$$
\xymatrix{
0 \arrow[r] &
\CU_- \arrow[r] \arrow[d]^{\xi_i} &
V\otimes\CO_{W_i} \arrow[r] \arrow[d] &
\CU_-^* \arrow[r] \arrow[d] &
0 \\
0 \arrow[r] &
\CU_+^* \arrow[r] &
\CU_+^* \arrow[r] &
0 \arrow[r] &
0 \\
}
$$
Note that the middle vertical arrow is surjective 
and its kernel is $\CU_+$. Hence, the long exact sequence
of kernels and cokernels gives
$0 \to \Ker \xi_i \to \CU_+ \to \CU_-^* \to \Coker \xi_i \to 0$.
Moreover, it is clear that the map $\CU_+ \to \CU_-^*$ in this
sequenece coincides with the dual map $\xi_i^*$. It follows immediately that 
$\Ker \xi_i \cong \Ker \xi_i^*$. On the other hand, it is clear that 
$\Ker \xi_i^* \cong (\Coker \xi_i)^*$.
\end{proof}

Let $\CE_i$ denote the cokernel of $\xi_i$ on $W_i$. It follows that $\CE_i$ is 
a rank $2$ vector bundle on $W_i$ and we have an exact sequence
\begin{equation}\label{ce}
0 \to \CE_i^* \to \CU_- \stackrel{\xi_i}{\to} \CU_+^* \to \CE_i \to 0.
\end{equation}
Dualizing, we obtain another sequence
\begin{equation}\label{ce2}
0 \to \CE_i^* \to \CU_+ \stackrel{\xi_i^*}{\to} \CU_-^* \to \CE_i \to 0.
\end{equation}

\begin{lemma}\label{cherns}
The Chern classes of bundles $\CE_i$ are given by the following formulas
$$
\begin{array}{ll}
\sfc_1(\CE_1) = H_X + H_{C^\vee}, &
\sfc_2(\CE_1) = \frac7{12} H_X H_{C^\vee} + 5 L_X + \eta,
\smallskip\\
\sfc_1(\CE_2) = H_S + H_{S^\vee}, &
\sfc_2(\CE_2) = \frac7{12} H_S H_{S^\vee} + 5 P_S + 5 P_{S^\vee},
\end{array}
$$
with $\eta\in \big(H^3(X,\C)\otimes H^1(C^\vee,\C)\big) \cap 
H^4(X\times C^\vee,\ZZ)$. 
\end{lemma}
\begin{proof}
It follows from~(\ref{ce}) that 
$$
\sfch(\CU_+^*) - \sfch(\CU_-) = 
\sfch(\CE_i) - \sfch(\CE_i^*) = 
2\sfch_1(\CE_i) + 2\sfch_3(\CE_i).
$$
This allows to compute 
$$
\begin{array}{ll}
\sfch_1(\CE_1) = H_X + H_{C^\vee}, &
\sfch_3(\CE_1) = -\frac12 P_X,
\smallskip\\
\sfch_1(\CE_2) = H_S + H_{S^\vee}, &
\sfch_3(\CE_2) = -\frac12 P_S -\frac12 P_{S^\vee}.
\end{array}
$$
Further, it is clear that $\sfc_1(\CE_i) = \sfch_1(\CE_i)$, and 
by K\"unneth formula we have
$$
\sfc_2(\CE_1) = a_1 H_X H_{C^\vee} + b_1 L_X + \eta,\qquad
\sfc_2(\CE_2) = a_2 H_S H_{S^\vee} + b_2 P_S + c_2 P_{S^\vee},
$$ 
for some $a_1,b_1,a_2,b_2,c_2\in\QQ$, 
$\eta\in\big(H^3(X,\C)\otimes H^1(C^\vee,\C)\big) \cap H^4(X\times C^\vee,\ZZ)$.
Further, since the correspondence $S \leftrightarrow S^\vee$ 
is symmetric, it is clear that $c_2 = b_2$. Finally, $a_i$ and $b_i$ 
can be found from the equality 
$3\sfc_1(\CE_i)\sfc_2(\CE_i) = \sfch_1(\CE_i)^3 - 6\sfch_3(\CE_i)$.
\end{proof}

\begin{remark}
Using the Riemann--Roch formula on $X \times C^{\vee}$ 
one can compute $\eta^2 = 14$.
\end{remark}

\begin{corollary}\label{uf}
The bundle $\CE_1$ {\rm(}resp.\ $\CE_2${\rm)} is a universal family 
of rank $2$ vector bundles with $\sfc_1=H_X$, $\sfc_2=5L_X$ on $X$ 
{\rm(}resp.\ with $\sfc_1=H_S$, $\sfc_2=5P_S$ on $S${\rm)}.
\end{corollary}
\begin{proof}
For every $y\in C^\vee$ (resp.\ $y\in S^\vee$) we denote by $\CE_{1y}$
the fiber of $\CE_1$ over $X\times y$ and by $\CE_{2y}$ the fiber of 
$\CE_2$ over $S\times y$. It will be shown in lemmas~\ref{he} and \ref{hes} 
below that all bundles $\CE_{1y}$ on $X$ for $y\in C^\vee$ and all bundles 
$\CE_{2y}$ on $S$ for $y\in S^\vee$ are stable, hence there exist morphisms 
$$
f_1:C^\vee\to\CM_X(2,H_X,5L_X),\qquad
f_2:S^\vee\to\CM_S(2,H_S,5P_S),
$$ 
to the moduli spaces of rank $2$ vector bundles on $X$ and $S$
with the indicated rank and Chern classes, such that
\begin{equation}\label{fstar}
\CE_1 = (\id_X\times f_1)^*\CE'_1\otimes q_1^*\CL_1,\qquad
\CE_2 = (\id_S\times f_2)^*\CE'_2\otimes q_2^*\CL_2, 
\end{equation}
where $\CE'_1$ and $\CE'_2$ are universal families 
on $X\times\CM_X(2,H_X,5L_X)$ and $S\times\CM_S(2,H_S,5P_S)$ respectively, 
$q_1:X\times C^\vee\to C^\vee$ and $q_2:S\times S^\vee \to S^\vee$ 
are the projections, and $\CL_1$ and $\CL_2$ are line bundles 
on $C^\vee$ and $S^\vee$ respectively.

It is easy to see that the maps $f_1$ and $f_2$ coincide with 
the maps $\rho$ constructed in \cite{IM}, section~4. Hence they 
are isomorphisms, and the bundles $\CE_1$ and $\CE_2$ are universal.
\end{proof}

Let $\alpha:S\to X$ and $\beta:C^\vee\to S^\vee$ denote the embeddings and put 
$\lambda_1 = \alpha\times\id_{C^\vee}$,
$\lambda_2 = \id_{S}\times\beta$,
$\mu_1 = \id_{X}\times\beta$,
$\mu_2 = \alpha\times\id_{S^\vee}$,
$\nu = \alpha\times\beta$.
Then we have a commutative diagram
$$
\xymatrix{
& S \times C^\vee \ar[dl]_{\lambda_1} \ar[dr]^{\lambda_2} \ar[dd]^{\nu}\\
X \times C^\vee \ar[dr]^{\mu_1} && S \times S^\vee \ar[dl]_{\mu_2} \\
& X \times S^\vee
}
$$

\begin{lemma}\label{defe}
We have canonical isomorphism $\lambda_1^*\CE_1 = \lambda_2^*\CE_2$.
\end{lemma}
\begin{proof}
The claim is clear since $\lambda_i^*\CE_i$ is the cokernel
of $\lambda_i^*\xi_i$, and $\lambda_1^*\xi_1 = \lambda_2^* \xi_2$
by definition of $\xi_i$.
\end{proof}

We denote the bundle $\lambda_1^*\CE_1 = \lambda_2^*\CE_2$
on $S \times C^\vee$ by $\CE$.

Consider the product $\TW = X \times S^\vee$ and the composition
$$
\txi : \CU_i \to V \otimes \CO_\TW \cong V^*\otimes \CO_\TW \to \CU_+^*.
$$

It is clear that
\begin{equation}\label{tc}
\mu_i^*\txi = \xi_i.
\end{equation}

\begin{lemma}\label{rc5}
The rank of $\txi$ equals $5$ at 
$X\times S^\vee \setminus 
\big( \mu_1(X\times C^\vee) \cup \mu_2(S\times S^\vee) \big )$.
\end{lemma}
\begin{proof}
Follows from lemma~\ref{di}.
\end{proof}

Let $\TE$ denote the cokernel of $\txi$. 

\begin{lemma}\label{te}
We have exact sequences on $X\times S^\vee$
$$
0 \to \CU_- \stackrel{\txi}{\to} \CU_+^* \to \TE \to 0,
$$
$$
0 \to \TE \to {\mu_1}_*\CE_1 \oplus {\mu_2}_*\CE_2 \to \nu_* \CE \to 0.
$$
\end{lemma}
\begin{proof}
The first sequence is exact by lemma~\ref{rc5} and definition of $\TE$.
To verify exactness of the second sequence we note that $\mu_i^*\TE = \CE_i$
by~(\ref{tc}) and~(\ref{ce}), and the canonical surjective maps 
$\TE \to {\mu_i}_*\mu_i^*\TE ={\mu_i}_*\CE_i$ glue to a surjective map
$\TE \to \Ker({\mu_1}_*\CE_1 \oplus {\mu_2}_*\CE_2 \to \nu_*\CE)$.
On the other hand, it is easy to check that the Chern characters of
$\TE$ and $\Ker({\mu_1}_*\CE_1 \oplus {\mu_2}_*\CE_2 \to \nu_*\CE)$
coincide, hence
$\TE \cong \Ker({\mu_1}_*\CE_1 \oplus {\mu_2}_*\CE_2 \to \nu_*\CE)$
and we are done.
\end{proof}

\begin{corollary}\label{ce12}
We have exact sequence on $X\times S^\vee$
$$
0 \to {\mu_1}_*\CE_1 \otimes \CO(-H_X) \to \TE \to {\mu_2}_*\CE_2 \to 0.
$$
\end{corollary}

\section{Cohomological computations}

\begin{lemma}\label{hu}
The pair $(\CU_+,\CO_X)$ in $\D^b(X)$ is exceptional. In orther words,
$$
\begin{array}{l}
\Ext^k(\CU_+,\CU_+) = H^k(X,\CU_+^*\otimes\CU_+) = 
\Ext^k(\CO_X,\CO_X) = H^k(X,\CO_X) = 
\begin{cases} \C, & \text{for $k=0$}\\0 , & \text{for $k\ne 0$}\end{cases}\\
H^\bullet(X,\CU_+)=0.
\end{array}
$$
\end{lemma}
\begin{proof}
Recall that $X$ is a complete intersection 
$X = \PP(A_9)\cap\LGr_+(V) \subset \PP(\SPV)$ and $\SPV/A_9 = B_7^*$.
Hence $X\subset\LGr_+(V)$ is the zero locus of a section of the vector bundle 
$B_7^*\otimes\CO_{\LGr_+(V)}(H_{\LGr_+(V)})$. Therefore, the Koszul 
complex $\Lambda^\bullet(B_7^*\otimes\CO_{\LGr_+(V)}(H_{\LGr_+(V)}))$ 
is a resolution of the structure sheaf $\CO_X$ on $\LGr_+(V)$. In other 
words, we have an exact sequence
$$
0 \to \Lambda^7(B_7\otimes\CO_{\LGr_+(V)}(-H_{\LGr_+(V)})) \to \!\ldots\! \to
\Lambda^1(B_7\otimes\CO_{\LGr_+(V)}(-H_{\LGr_+(V)})) \to 
\CO_{\LGr_+(V)} \to \CO_X \to 0.
$$
Tensoring it by $\CU_+$ and $\CU_+^*\otimes\CU_+$ we see that
it suffices to compute
$H^\bullet(\LGr_+(V),F(-kH_{\LGr_+(V)}))$ for $F = \CO_{\LGr_+(V)}$,
$F = \CU_+$ and $F = \CU_+^*\otimes\CU_+$ and $0\le k\le 7$.
These cohomologies are computed by Borel--Bott--Weil Theorem \cite{D},
since all the bundles under the question are the pushforwards 
of equivariant line bundles on the flag variety 
of the spinor group $\Spin(V)$.
\end{proof}

Since the canonical class of $X$ equals $-H_X$, the Serre duality on $X$ gives

\begin{corollary}\label{hus}
We have 
$$
H^\bullet(X,\CU_+^*(-H_X))=0,\qquad
H^k(X,\CU_+\otimes\CU^*_+(-H_X)) = 
\begin{cases} \C, & \text{for $k=3$}\\0 , & \text{for $k\ne 3$}\end{cases}.
$$
\end{corollary}

\begin{lemma}\label{he}
For any $y\in C^\vee$ we have 
$H^p(X,{\CE_1}_y(-H_X)) = H^p(X,{\CE_1}_y\otimes\CU_+^*(-H_X)) = 0$.
In particular, $\CE_{1y}$ is stable.
\end{lemma}
\begin{proof}
Recall that by definition $C^\vee = \LGr_-(V)\cap\PP(B_7)$, and
$X = \LGr_+(V)\cap\PP(A_9)$ with $A_9 = B_7^\perp$. 
Choose a hyperplane $\PP(B_6)\subset \PP(B_7)$  such that
$\PP(B_6)$ intersects $C^\vee$ transversally and doesn't contain $y$.
Take $A_{10} = B_6^\perp$ and consider $\HX = \LGr_+(V)\cap\PP(A_{10})$. 
Then the arguments of lemma~\ref{smoothness} show that $\HX$ is a smooth 
Fano fourfold of index 2 containing $X$ as a hyperplane section. 
Moreover, the arguments similar to that of lemma~\ref{te} show that
the composition of morphisms on $\HX$
$$
\hxi:{\CU_-}_y\otimes\CO_\HX \to V\otimes\CO_\HX \to \CU_+^*
$$
is injective and its cokernel is isomorphic to the pushforward of
${\CE_1}_y$ via the embedding $i:X\to\HX$. In other words, we have
the following exact sequence on $\HX$:
\begin{equation}\label{x1}
0 \to {\CU_-}_y\otimes\CO_\HX \to \CU_+^* \to i_* {\CE_1}_y \to 0,
\end{equation}
On the other hand, using Borel--Bott--Weil Theorem and
the Koszul resolution of $\HX\subset\LGr_+(V)$ along the lines of 
lemma~\ref{hu} one can compute
$$
H^\bullet(\HX,\CU_+^*\otimes\CU_+^*(-H_\HX)) = 
H^\bullet(\HX,\CU_+^*(-H_\HX)) = 
H^\bullet(\HX,\CO_\HX(-H_\HX)) = 0
$$
and the claim follows from the cohomology sequences of~(\ref{x1}) 
twisted by $\CO_\HX(-H_\HX)$ and $\CU_+^*(-H_\HX)$ respectively,
since 
$$
H^\bullet(X,{\CE_1}_y(-H_X)) = H^\bullet
(\HX,i_*{\CE_1}_y(-H_\HX)),
\quad
H^\bullet(X,{\CE_1}_y\otimes\CU_+^*(-H_X)) = 
H^\bullet(\HX,i_*{\CE_1}_y\otimes\CU_+^*(-H_\HX)).
$$
\end{proof}

\begin{lemma}\label{he2h}
For any $y\in C^\vee$ we have $H^1(X,\CE_{1y}(-2H_X)) = 0$.
\end{lemma}
\begin{proof}
Restricting exact sequence~$(\ref{ce})$ to 
$X = X\times \{y\} \subset X\times C^\vee$,
twisting it by $\CO_X(-H_X)$ and taking into account lemma~\ref{cherns}
we obtain exact sequence
$$
0 \to \CE_{1y}(-2H_X) \to {\CU_-}_y\otimes\CO_X(-H_X) \to 
\CU_+^*(-H_X) \to \CE_{1y}(-H_X) \to 0.
$$
It follows from corollary~\ref{hus} and lemma~\ref{he} that
$H^1(X,\CE_{1y}(-2H_X)) = {\CU_-}_y\otimes H^1(X,\CO_X(-H_X))$,
but using Serre duality we have $H^1(X,\CO_X(-H_X)) = H^2(X,\CO_X)^* = 0$
by lemma~\ref{hu}.
\end{proof}

\begin{lemma}\label{hes}
For any $y\in S^\vee$ we have $H^0(S,\CE_{2y}(-H_S))=0$.
In particular, $\CE_{2y}$ is stable.
\end{lemma}
\begin{proof}
For $y\in C^\vee$ we have $\CE_{2y} = {\CE_{1y}}_{|S}$, hence 
the claim follows from exact sequence
$$
H^0(X,\CE_{1y}(-H_X)) \to 
H^0(S,\CE_{2y}(-H_S)) \to 
H^1(X,\CE_{1y}(-2H_X)),
$$
since the first term vanishes by lemma~\ref{he},
and the third term vanishes by lemma~\ref{he2h}.

Now we note that while $S$ (and hence $S^\vee$) is fixed we can take
for $C^\vee$ any smooth hyperplane section of $S^\vee$, consider
the corresponding smooth $X \supset S$, and repeat the above arguments
in this situation. Since any point $y\in S^\vee$ lies on a smooth 
hyperplane section, these arguments prove the claim for all $y\in S^\vee$.
\end{proof}

\begin{corollary}\label{heu}
For any $y\in C^\vee$ we have 
$H^\bullet(X,{\CE_1}_y\otimes\CU_+(-H_X)) = 0$.
\end{corollary}
\begin{proof}
Tensor exact sequence $0 \to \CU_+ \to V\otimes\CO_X \to \CU_+^* \to 0$
with ${\CE_1}_y(-H_X)$ and consider the cohomology sequence.
\end{proof}

\section{Derived categories}

Consider the kernel functors taking $\CE_1$ and $\CE_2$ for kernels:
$$
\Phi_1:\D^b(C^\vee) \to \D^b(X),\qquad
\Phi_2:\D^b(S^\vee) \to \D^b(S),\qquad
\Phi_i(-) = R{p_i}_*(L{q_i}^*(-)\otimes\CE_i),
$$
where $p_i$ and $q_i$ are the projections onto 
the first and the second factors:
$$
\xymatrix{
& 
X \times C^\vee \ar[dl]_{p_1} \ar[dr]^{q_1} 
&&&&
S \times S^\vee \ar[dl]_{p_2} \ar[dr]^{q_2} \\
X && C^\vee && S && S^\vee
}
$$

\begin{theorem}
The functor $\Phi_i$ is fully faithful.
\end{theorem}
\begin{proof}
According to the result of Bondal and Orlov \cite{BO} it suffices to check
that for the structure sheaves of any two points $y_1,y_2\in C^\vee$
(resp. $y_1,y_2\in S^\vee$) and all $p\in\ZZ$ we have
$$
\Ext^p(\Phi_i(\CO_{y_1}),\Phi_i(\CO_{y_2})) =
\Ext^p(\CO_{y_1},\CO_{y_2}).
$$
But clearly $\Phi_i(\CO_{y_k}) = \CE_{iy_k}$
and it remains to apply corollary~\ref{uf} and theorem~\ref{im2}.
\end{proof}

\begin{corollary}\label{phi2}
The functor $\Phi_2:\D^b(S^\vee) \to \D^b(S)$ is an equivalence.
\end{corollary}
\begin{proof}
Any fully faithful functor between the derived categories 
of K3 surfaces is an equivalence, see \cite{Br1}.
\end{proof}

Consider the following diagram
$$
\xymatrix{
X \ar@{.}[r]^{\CE_1} &
C^\vee \ar[d]^{\beta} \\
S \ar[u]^{\alpha} \ar@{.}[r]^{\CE_2} &
S^\vee 
}
$$
where the dotted line connecting two varieties means that we consider 
the corresponding kernel on their product. This diagram induces a diagram
of functors
$$
\xymatrix{
\D^b(X) \ar[d]_{\alpha^*} &
\D^b(C^\vee) \ar[d]^{\beta_*} \ar[l]_{\Phi_1} \\
\D^b(S) &
\D^b(S^\vee) \ar[l]_{\Phi_2}
}
$$
which is commutative by lemma~\ref{defe}, since the functor 
$\alpha^*\circ\Phi_1$ is given by the kernel $\lambda_1^*\CE_1$,
and the functor $\Phi_2\circ\beta_*$ is given by the kernel 
$\lambda_2^*\CE_2$.

Let $\Phi_1^*:\D^b(X) \to \D^b(C^\vee)$ and 
$\Phi_2^*:\D^b(S) \to \D^b(S^\vee)$ denote 
the left adjoint functors. The standard computation
shows that these functors are given by the kernels
$$
\CE_1^*(-H_X)[3] = \CE_1(-2H_X-H_{C^\vee})[3]\quad\text{on $X\times C^\vee$,} 
\quad\text{and}\quad
\CE_2^*[2] = \CE_2(-H_S-H_{S^\vee})[2]\quad\text{on $S\times S^\vee$}
$$ 
respectively. Consider the following diagram
$$
\xymatrix{
\D^b(X) \ar[d]_{\alpha^*} \ar[r]^{\Phi_1^*} &
\D^b(C^\vee) \ar[d]^{\beta_*} \\
\D^b(S) \ar[r]^{\Phi_1^*} &
\D^b(S^\vee) 
}
$$
This diagram is no longer commutative, however, the following proposiotion
shows that it becomes commutative if one replaces $\D^b(X)$ by its 
subcategory ${}^\perp\langle\CU_+,\CO_X\rangle \subset \D^b(X)$.

\begin{proposition}\label{commdiag}
The functors $\beta_*\circ\Phi_1^*$ and
$\Phi_2^*\circ\alpha^*:\D^b(X) \to \D^b(S^\vee)$
are isomorphic on the subcategory 
${}^\perp\langle\CU_+,\CO_X\rangle \subset \D^b(X)$.
\end{proposition}
\begin{proof}
It is clear that the functors $\beta_*\circ\Phi_1^*$ and
$\Phi_2^*\circ\alpha^*$ are given by the kernels 
${\mu_1}_*\CE_1(-2H_X-H_{S^\vee})[3]$ 
and ${\mu_2}_*\CE_2(-H_X-H_{S^\vee})[2]$ on $X\times S^\vee$
respectively. Considering the helix of the exact sequence of 
corollary~\ref{ce12} twisted by $\CO(-H_X-H_{S^\vee})$ we see 
that there exists a distinguished triangle
$$
{\mu_2}_*\CE_2(-H_X-H_{S^\vee})[2] \to 
{\mu_1}_*\CE_1(-2H_X-H_{S^\vee})[3] \to
\TE(-H_X-H_{S^\vee})[3].
$$ 
It remains to show that a kernel functor $\D^b(X)\to \D^b(S^\vee)$ 
given by the kernel $\TE(-H_X-H_{S^\vee})$ vanishes on the triangulated 
subcategory ${}^\perp\langle\CU_+,\CO_X\rangle \subset \D^b(X)$.

Note that lemma~\ref{te} implies that $\TE(-H_X-H_{S^\vee})$ 
is isomorphic to a cone of the morphism 
$\txi:\CU_-(-H_X-H_{S^\vee}) \to \CU_+^*(-H_X-H_{S^\vee})$
on $X\times S^\vee$, so it suffices to check that the kernel
functors $\D^b(X) \to \D^b(S^\vee)$ given by the kernels
$\CU_-(-H_X-H_{S^\vee})$ and $\CU_+^*(-H_X-H_{S^\vee})$ on $X\times S^\vee$
vanish on the triangulated subcategory 
${}^\perp\langle\CU_+,\CO_X\rangle \subset \D^b(X)$.
Let $\tp:X\times S^\vee \to X$ and $\tq:X\times S^\vee \to S^\vee$
denote the projections. The straightforward computation 
using the projection formula and the Serre duality on $X$ shows that
for any object $F\in \D^b(X)$ we have
\begin{multline*}
\Phi_{\CU_-(-H_X-H_{S^\vee})}(F) =
R\tq_*(L\tp^*(F) \otimes \CU_-(-H_X-H_{S^\vee})) = \\ =
\RH(X,F(-H_X)) \otimes \CU_-(-H_{S^\vee}) =
\RHom(F,\CO_X)^*\otimes \CU_-(-H_{S^\vee}),
\end{multline*}
\begin{multline*}
\Phi_{\CU_+^*(-H_X-H_{S^\vee})}(F) = 
R\tq_*(L\tp^*(F) \otimes \CU_+^*(-H_X-H_{S^\vee})) = \\ =
\RH(X,F\otimes \CU_+^*(-H_X)) \otimes \CO_{S^\vee}(-H_{S^\vee}) = 
\RHom(F,\CU_+)^*\otimes \CO_{S^\vee}(-H_{S^\vee}).
\end{multline*}
In particular, the above kernel functors vanish for all objects 
$F\in {}^\perp\langle\CU_+,\CO_X\rangle \subset \D^b(X)$ and we are done.
\end{proof}

\begin{theorem}
We have a semiorthogonal decomposition
\begin{equation}\label{decomp}
\D^b(X) = \big\langle\CU_+,\CO_X,\Phi_1(\D^b(C^\vee))\big\rangle.
\end{equation}
\end{theorem}
\begin{proof}
It is clear that $\CO_X$ is an exceptional bundle, and
$\CU_+$ is an exceptional bundle by lemma~\ref{hu}.
Now, let us verify the semiorthogonality. Indeed, 
$$
\Ext^\bullet(\CO_X,\CU_+) = H^\bullet(X,\CU_+) = 0
$$
by lemma~\ref{hu}. Moreover, denoting by $\Phi_1^!:\D^b(X)\to\D^b(C^\vee)$
the right adjoint functor and taking either $F = \CO_X$, or $F = \CU_+$ 
we see that
$$
\Ext^\bullet(\CO_y,\Phi_1^!(F)) = 
\Ext^\bullet(\Phi_1(\CO_y),F) =
\Ext^\bullet(\CE_y,F) = 
H^p(X,\CE_y^*\otimes F) = 
H^p(X,\CE_y\otimes F(-H_X)) = 0
$$
by lemma~\ref{he} and corollary~\ref{heu}. Hence $\Phi_1^!(F) = 0$, 
since $\{\CO_y\}_{y\in C^\vee}$ is a spanning class (see~\cite{Br1})
in $\D^b(C^\vee)$, hence
$$
\Ext^\bullet(\Phi_1(G),F) =
\Ext^\bullet(G,\Phi_1^!(F)) = 
\Ext^\bullet(G,0) = 0
$$
for all $G\in\D^b(C^\vee)$.

It remains to check that $\D^b(X)$ is generated by $\CU_+$, $\CO_X$, 
and $\Phi_1(\D^b(C^\vee))$ as a triangulated category. Indeed, assume that
$F \in {}^\perp\langle\CU_+,\CO_X,\Phi_1(\D^b(C^\vee))\rangle$.
Since $F \in {}^\perp\Phi_1(\D^b(C^\vee))$ we have $\Phi_1^*(F) = 0$.
On the other hand, since $F \in {}^\perp\langle\CU_+,\CO_X\rangle$ 
we have by proposition~\ref{commdiag}
$$
\Phi_2^*\circ\alpha^*(F) = \beta_*\circ\Phi_1^*(F) = 0.
$$
But $\Phi_2$ is an equivalence by corollary~\ref{phi2},
hence $\Phi_2^*$ is an equivalence, hence $\alpha^*(F) = 0$.

Now we note, that while $X$ (and hence $C^\vee$) is fixed, we can take 
for $S$ any smooth hyperplane section of $X$. Then the above arguments 
imply that for any 
$F \in {}^\perp\langle\CU_+,\CO_X,\Phi_1(\D^b(C^\vee))\rangle \subset \D^b(X)$
its restriction to any smooth hyperplane section is isomorphic to zero.
Thus the proof is finished by the following lemma.
\end{proof}

\begin{lemma}
If $X$ is a smooth algebraic variety and $F$ is a complex
of coherent sheaves on $X$ which restriction to every smooth
hyperplane section of $X$ is acyclic, then $F$ is acyclic.
\end{lemma}
\begin{proof}
Assume that $F$ is not acyclic and let $k$ be the maximal integer 
such that $\CH^k(F) \ne 0$. Let $x\in X$ be a point in the support 
of the sheaf $\CH^k(F)$. Choose a smooth hyperplane section $j:S\subset X$
passing through $x$. Since the restriction functor $j^*$ is right-exact
it is clear that $\CH^k(Lj^*F) \ne 0$, a contradiction.
\end{proof}

\section{Application: the Fano surface of conics}

Let $F_X$ denote the Fano surface of conics (rational curves 
of degree $2$) on $X$.

\begin{lemma}\label{upr}
If $R\subset X$ is a conic then 
${\CU_+}_{|R} \cong \CO_R \oplus \CO_R(-1)^{\oplus 4}$.
\end{lemma}
\begin{proof}
Since $\CU_+$ is a subbundle of the trivial vector bundle $V\otimes\CO_X$,
and since we have $r({\CU_+}_{|R}) = 5$, $\deg({\CU_+}_{|R}) = -4$
we have ${\CU_+}_{|R} \cong \oplus_{j=1}^5 \CO_R(-u_j)$, where
$u_j\ge 0$ and $\sum u_j = 4$. Thus it suffices to check that
$\dim H^0(R,{\CU_+}_{|R}) = 1$. Actually, $\dim H^0(R,{\CU_+}_{|R}) \ge 1$
follows from above, so it remains to show that
$\dim H^0(R,{\CU_+}_{|R}) \ge 2$ is impossible. 

Indeed, assume $\dim H^0(R,{\CU_+}_{|R}) \ge 2$. Choose a $2$-dimensional
subspace $U \subset H^0(R,{\CU_+}_{|R}) \subset H^0(R,V\otimes\CO_R) = V$
and consider $V' = U^\perp/U$. Then $\LGr_+(V') \subset \LGr_+(V)$,
and it is clear that $R \subset X' := \LGr_+(V') \cap X$.
Since $X$ is a plane section of $\LGr_+(V)$, therefore
$X'$ is a plane section of $\LGr_+(V')$. But $V' = \C^6$, hence
$\LGr_+(V') \cong \PP^3$. But a plane section of $\PP^3$ containing
a conic contains a plane $\PP^2$, hence $X'$ contains $\PP^2$,
hence $X$ contains $\PP^2$ which contradicts Lefschetz theorem for $X$.
\end{proof}

\begin{lemma}\label{capum}
We have $\bigcap_{y\in C^\vee}{\CU_-}_y = 0$.
\end{lemma}
\begin{proof}
Assume that $0\ne v\subset \bigcap_{y\in C^\vee}{\CU_-}_y$ and
consider $V'' = v^\perp/\C v$. 
Then $C^\vee \subset \LGr_-(V'') \subset \LGr_-(V)$.
Moreover, since $C^\vee$ is a plane section of $\LGr_-(V)$,
hence $C^\vee$ is a plane section of $\LGr_-(V'')$.
Further, $V'' = \C^8$, hence $\LGr_-(V'')$ is a quadric, 
and a curve which is a plane section of a quadric 
is a line or a conic. But $C^\vee$ is neither.
\end{proof}

\begin{theorem}
We have $F_X \cong S^2C^\vee$.
\end{theorem}
\begin{proof}
Let $R\subset X$ be a conic and consider a decomposition 
of its structure sheaf with respect to the semiorthogonal
decomposition
$$
\D^b(X) = \langle\CO_X,\CU_+^*,\Phi_1(\D^b(C^\vee))\rangle,
$$
obtained from the decomposition~(\ref{decomp}) by mutating $\CU_+$
through $\CO_X$. To this end we compute
$$
\Ext^p(\CO_R,\CO_X) = 
H^{3-p}(X,\omega_X\otimes\CO_R)^* = 
H^{3-p}(R,\omega_R)^* = 
\begin{cases}
\C, & \text{if $p=2$}\\
0, & \text{otherwise}
\end{cases}
$$
\begin{multline*}
\Ext^p(\CO_R,\CU_+) = 
\Ext^{3-p}(\CU_+,\omega_X\otimes\CO_R)^* = 
\Ext^{3-p}({\CU_+}_{|R},\omega_R)^* = \\ =
H^{3-p}(R,\CU_+^*\otimes\omega_R)^* = 
H^{p-2}(R,{\CU_+}_{|R}) =
\begin{cases}
\C, & \text{if $p=2$}\\
0, & \text{otherwise}
\end{cases}
\end{multline*}
by lemma~\ref{upr}.
Hence the decomposition gives the following exact sequence
\begin{equation}\label{decc}
0 \to \CO_X \to \CU_+^* \to \Phi_1(\Phi_1^!(\CO_R)) \to \CO_R \to 0,
\end{equation}
where 
$$
\Phi_1^!:\D^b(X) \to \D^b(C^\vee),\qquad
\Phi_1^!(-) = R{q_1}_*(Lp_1^*(-)\otimes\CE_1^*(H_{C^\vee}))[1],
$$
is the right adjoint to $\Phi_1$ functor.

\begin{lemma}
$\Phi_1^!(\CO_R)$ is a pure sheaf.
\end{lemma}
\begin{proof}
In order to understand 
$\Phi_1^!(\CO_R) = R{q_1}_*(Lp_1^*(\CO_R)\otimes\CE_1^*(H_{C^\vee}))[1] 
\in \D^b(C^\vee)$,
we investigate
$H^\bullet(X,{\CE_1^*}_y\otimes\CO_R) = H^\bullet(R,{{\CE_1^*}_y}_{|R})$
for all $y\in C^\vee$. The sheaf ${\CE_1^*}_y$ by~(\ref{ce}) is a subsheaf 
of the trivial vector bundle ${\CU_-}_y\otimes\CO_X$, therefore
$H^0(R,{{\CE_1^*}_y}_{|R}) \subset {\CU_-}_y \subset V$.
On the other hand, by~(\ref{ce2}) we have ${\CE_1^*}_y={\CE_1}_y(-H_X)$, 
is a subsheaf of $\CU_+$, hence
$H^0(R,{{\CE_1^*}_y}_{|R}) \subset H^0(R,{\CU_+}_{|R}) = \C \subset V$.
Therefore, if $H^0(R,{{\CE_1^*}_y}_{|R})\ne 0$ for all $y\in C^\vee$
then $\bigcap_{y\in C^\vee}{\CU_-}_y \ne 0$ which is false 
by lemma~\ref{capum}. Thus for generic $y\in C^\vee$ we have 
$H^0(R,{{\CE_1^*}_y}_{|R}) = 0$, hence
$R^0{q_1}_*(Lp_1^*(\CO_R)\otimes\CE_1^*(H_{C^\vee})) = 0$.
On the other hand, since $R$ is $1$-dimensional we have
$R^k{q_1}_*(Lp_1^*(\CO_R)\otimes\CE_1^*(H_{C^\vee})) = 0$ for $k\ne 0,1$.
Hence $\Phi_1^!(\CO_R)$ is a pure sheaf.
\end{proof}

\begin{corollary}
$\Phi_1^!(\CO_R)$ is an artinian sheaf of length $2$ on $C^\vee$.
\end{corollary}
\begin{proof}
Computation of the Chern character of $\Phi_1^!(\CO_R)$
via the Grothendieck--Riemann--Roch.
\end{proof}

It follows from above that $\Phi_1^!(\CO_R)$ is either 
the structure sheaf of a length $2$ subscheme in $C^\vee$, or
$\Phi_1^!(\CO_R) = \CO_y \oplus \CO_y$ for some 
$y\in C^\vee$. We claim that the second never happens.
To this end we need the following
                                            
\begin{lemma}\label{phistaru}
We have $\Phi_1^*(\CU_+^*) = \CO_{C^\vee}$.
\end{lemma}
\begin{proof}
It is clear that
$$
\Phi_1^*(\CU_+^*) = R{q_1}_*(L{p_1}^*(\CU_+^*)\otimes\CE_1^*(-H_X))[3] =
R{q_1}_*(\CE_1\otimes\CU_+^*(-2H_X-H_{C^\vee}))[3].
$$
On the other hand, tensoring (\ref{ce2}) with $\CU_+^*(-H_X)$
we obtain exact sequence
$$
0 \to 
\CE_1\otimes\CU_+^*(-2H_X-H_C^\vee) \to 
\CU_+\otimes\CU_+^*(-H_X) \to
\CU_-^*\otimes\CU_+^*(-H_X) \to
\CE_1\otimes\CU_+^*(-H_X) \to 0.
$$
Lemma~\ref{he} implies that 
$R^\bullet{q_1}_*(\CE_1\otimes\CU_+^*(-H_X)) = 0$ 
and corollary~\ref{hus} implies that
$$
R^\bullet{q_1}_*(\CU_-^*\otimes\CU_+^*(-H_X)) = 0,\qquad
R^k{q_1}_*(\CU_+\otimes\CU_+^*(-H_X)) = 
\begin{cases}
\CO_{C^\vee}, & \text{for $k=3$}\\
0, & \text{for $k\ne 3$}
\end{cases}
$$
and the claim follows from the spectral sequence.
\end{proof}

\begin{lemma}
$\Phi_1^!(\CO_R) \ne \CO_y \oplus \CO_y$.
\end{lemma}
\begin{proof}
If the above would be true then the decomposition~(\ref{decc}) 
would take form
$$
0 \to \CO_X \to \CU_+^* \to {\CE_1}_y \oplus{\CE_1}_y \to \CO_R \to 0.
$$
On the other hand, it follows from lemma~\ref{phistaru} that
$$
\Hom(\CU_+^*,{\CE_1}_y) =
\Hom(\CU_+^*,\Phi_1(\CO_y)) = 
\Hom(\Phi_1^*(\CU_+^*),\CO_y) =
\Hom(\CO_{C^\vee},\CO_y) = \C,
$$
hence the map $\CU_+^* \to {\CE_1}_y \oplus {\CE_1}_y$ must have rank~$2$
and the above sequence is impossible.
\end{proof}

\begin{corollary}
$\Phi_1^!(\CO_R)$ is the structure sheaf 
of a length $2$ subscheme in $C^\vee$.
\end{corollary}

Thus the functor $\Phi_1^!$ induces a map $F_X \to S^2 C^\vee$.

Vice versa, if $Z$ is a length 2 subscheme in $C^\vee$ then
$$
\Hom(\CO_{C^\vee},\CO_Z) =
\Hom(\Phi_1^*(\CU_+^*),\CO_Z) =
\Hom(\CU_+^*,\Phi_1(\CO_Z)). 
$$
Therefore, the canonical projection $\CO_{C^\vee} \to \CO_Z$ induces
canonical morphism $f:\CU_+^*\to\Phi_1(\CO_Z)$. Its kernel, 
being a rank $1$ reflexive sheaf with $\sfc_1 = 0$,
must be isomorphic to $\CO_X$, and it is easy to show that its
cokernel is the structure sheaf of a conic.
Therefore, the map $\Phi_1^!:F_X \to S^2 C^\vee$ is an isomorphism.
\end{proof}

\end{document}